\documentclass[11pt]{amsart}
\usepackage{amsmath}
\newtheorem{thm}{Theorem}[section]
\newtheorem{prop}{Proposition}[section]
\newtheorem{lemma}{Lemma}[section]

\theoremstyle{remark}
\newtheorem*{remark}{Remarks}
\theoremstyle{definition}
\newtheorem*{definition}{Example}
\newcommand{\N}{{\mathbb N}}
\newcommand{\R}{{\mathbb R}}

\newcommand{\Z}{{\mathbb Z}}

\title[Sharp Bounds for Eigenvalues]{Sharp Bounds for Eigenvalues and Multiplicities
 on Surfaces of Revolution}
\author{Martin Engman} 
\begin{document}
\address{Department of Mathematics, University of New Mexico, Albuquerque, N.M. 87131}
\email{engman@math.unm.edu}
\begin{abstract}
{For surfaces of revolution diffeomorphic to $S^2$, it is proved that $(S^2,can)$ provides sharp upper bounds for the multiplicities of all of the distinct eigenvalues.
 We also find sharp upper bounds for all the distinct eigenvalues and show that an infinite sequence of these eigenvalues are bounded above by those of $(S^2,can)$.
 An example of such bounds for a metric with some negative curvature is presented.  }
\end{abstract}
\maketitle
\section{Introduction}

Upper bounds for the multiplicities of eigenvalues have been found by Cheng and Besson in \cite{ch1} and \cite{bes1}, respectively.
 Besson obtained the upper bound $m_k(g)\leq 4p+2k+1$ for the multiplicity of the $k$th eigenvalue of any compact Riemannian surface of genus $p$.
 Since these multiplicity bounds are for eigenvalues which are not necessarily distinct, even the bound $m_k(g)\leq 2k+1$ for $p=0$ fails to be sharp after $k=1$.
If, for $p=0$, one could obtain the same formula, but for the \underline{distinct} eigenvalues, then the bound would be sharp since the multiplicity of the $k$th distinct eigenvalue on $(S^2,can)$ is $2k+1$.
 However, without additional restrictions on the metric $g$, this problem seems to be quite difficult.
 We will restrict our attention to surface of revolution metrics in order to obtain a special case of this result via the methods of the author's previous article \cite{eng}.
Although one might still hope to obtain such a result for arbitrary metrics on $S^2$, there can be no generalization to higher dimensional spheres because of the results of \cite{bb} and \cite{ur}.

Another problem of considerable interest is that of bounding the eigenvalues themselves.
 Cheng (\cite{ch2}) obtained upper bounds and later, Li and Yau (\cite{ly}) obtained both upper and lower bounds in a very general setting.
 Again, these bounds are for the sequence of eigenvalues counted with their multiplicities and, therefore, may not be sharp.
 In this paper we obtain, for the special case of surfaces of revolution diffeomorphic to $S^2$, sharp upper bounds for the distinct eigenvalues.
 It is also shown that the upper bounds are achieved only for the constant curvature metric.
In general, these upper bounds are not, themselves, bounded by the eigenvalues of $(S^2,can)$.
 But it is shown that there exists an infinite subsequence of distinct eigenvalues which are bounded by the corresponding eigenvalues of $(S^2,can)$. 

The paper concludes with an example of a metric with some negative curvature whose spectrum is strictly less than that of the standard sphere and, in fact, diverges from that of the standard sphere.   

\section{Previous Results}

 In this section we present, without proof, those results from \cite{eng} which are necessary for our purposes here.

Let $(M,g)$ be a surface of revolution of surface area $4\pi$ which is diffeomorphic to the sphere.
 In \cite{eng} it was shown that on a certain chart the metric can take the form:
\begin{equation}
 g=\frac{1}{f(x)}dx \otimes dx + f(x)d\theta \otimes d\theta   \label{eq:met}
\end{equation}
where $(x,\theta) \in(-1,1) \times [0,2\pi)$ and $f(x)
 \in C^{\infty}(-1,1) \cap C^{0}[-1,1]$ satisfies the following conditions:
\begin{equation}
   f(x)>0  \; \; \forall x \in (-1,1),\;\;
   f(1)=f(-1)=0 \; \mbox{ and, }
   f'(-1)=-f'(1)=2.      \label{eq:f}
\end{equation}
This metric has Gauss curvature given by $K(
x) =(-1/2)f^{''}(x)$. The canonical (i.e. constant curvature) metric is obtained by taking $f(x)=1-x^{2}$ and the metric in this case is denoted by $can$.

Let $\Delta$ denote the Laplacian on $(M,g)$ and let $\lambda$ be any eigenvalue of $-\Delta$. 
We use the symbols $E_\lambda$, and $dimE_\lambda$ to mean the eigenspace for $\lambda$ and its multiplicity, respectively. In this paper $\lambda_m$ will always mean the $m$th \underline{distinct} eigenvalue. 
Since $S^1$ (parameterized here by $0 \leq \theta < 2\pi$) acts on $(M,g)$ by isometries, the orthogonal decomposition of $E_\lambda$ has the special form:

\begin{equation}
E_\lambda= \left\{ \begin{array}{ll}
 \; \;  \bigoplus_{j=1}^{l} ( e^{-ik_{j}\theta} V_{k_{j}}\oplus e^{ik_{j}\theta} V_{k_{j}}) & \mbox{if $dimE_{\lambda}$ is even} \\
 \left\{ \bigoplus_{j=1}^{l} ( e^{-ik_{j}\theta} V_{k_{j}}\oplus e^{ik_{j}\theta} V_{k_{j}}) \right\} \oplus V_{0} & \mbox{if $dimE_{\lambda}$ is odd}
                      \end{array} \right.  \label{eq:1} 
\end{equation}
where $k_{j} \in \N$ and $V_{k_{j}}$ is the one dimensional eigenspace of the ordinary differential operator 
\begin{equation}
L_{k_{j}}=-\frac{d}{dx}\left(f(x)\frac{d}{dx}\right) + \frac{k_{j}^{2}}{f(x)}.       \label{eq:odo}
\end{equation}
acting on functions which vanish at $1$ and $-1$.
Consequently, the spectrum of $-\Delta$ can be studied via the spectra $SpecL_{k}=\{\lambda_{k}^{1} < \lambda_{k}^{2} < \cdots <\lambda_{k}^{j} < \cdots \}$ for all $k \in \Z$.
 When $k=0$ the spectrum is called {\em $S^{1}$ invariant}, otherwise it is called {\em $k$-equivariant}.
 Each $L_{k}$ has a Green's operator, $\Gamma_{k}$, whose spectrum is, of course, $\{ 1/\lambda_{k}^{j} \}_{j=1}^{\infty}$ and whose trace is given by $tr \Gamma_{k}=\sum 1/\lambda_{k}^{j}$.
 Recall that for any constant curvature metric, $can$, on $S^{2}$, $dimE_{\lambda_{m}(can)}=2m+1$.  

\begin{prop}[See \cite{eng}] The following hold on any $(M,g)$ where $g$ is a surface of revolution metric of area $4\pi$, and $M$ is diffeomorphic to $S^{2}$.
  \begin{description}
   \item[i] ${\displaystyle Spec(-\Delta)=\cup_{k\in \Z} \left\{ SpecL_{k} \right\}.}$
   \item[ii] If $0<k<l$ then $\lambda_{k}^{1}<\lambda_{l}^{1}$. (Monotonicity) 
   \item[iii] ${\displaystyle tr \Gamma_{k} = \tfrac{1}{|k|}}$, for $k \neq 0$. (Trace Formula)
   \item[iv] Let $\lambda_{m}$ be the $m$th distinct eigenvalue of $-\Delta$.
Then $dimE_{\lambda_{m}}=2m+1$ for all $m$ if and only if $(M,g)$ is isometric to a sphere of constant curvature.
   \item[v] If $h$ and $h^2$ are both eigenfunctions of $-\Delta$, then the metric is of constant curvature $\equiv 1$.
  \end{description}
\end{prop}

\begin{remark} 
\begin{enumerate}
  \item The explicit form of the trace formula is:
   \begin{equation}
    \sum_{j=0}^{\infty} 1/\lambda_{k}^{j+1}=\frac{1}{|k|}. \label{eq:2}
   \end{equation}
This formula is the main ingredient in the proof of iv and will play an important role in what follows.
   \item Item iv is one of the main results of \cite{eng}.
 Since its publication, S.Y. Cheng has made some progress toward a more general result.
 He has proved that any metric on $S^{2}$ with $dimE_{\lambda_{m}}=2m+1$ is a Z\"{o}ll metric. 
\end{enumerate}
\end{remark}
\section{Sharp Bounds for the Multiplicities}

Before presenting the main theorems we prove a lemma which shows the relationship between the distinct eigenvalues of $-\Delta$ and those of the operators $L_k$.
\begin{lemma} We will assume the same hypotheses as Proposition 2.1. For all $k\geq 1$ and $j\geq 0$, $\lambda_{k+j} \leq \lambda^{j+1}_{k}$.
\end{lemma}
\begin{proof} Because of Proposition 2.1, i., ii., and the simplicity of the spectrum of $L_k$, there is a strictly increasing subsequence of eigenvalues of $-\Delta$, $\lambda_{1}^{1} < \lambda_{2}^{1} < \cdots <\lambda_{k}^{1} <\lambda_{k}^{2} < \cdots <\lambda_{k}^{j+1}$.
 As there are $(k+j)$ distinct eigenvalues in this list, we must have that $\lambda_{k}^{j+1}$ is at least the $(k+j)$th distinct eigenvalue of $-\Delta$.
 In other words,  $\lambda_{k+j} \leq \lambda^{j+1}_{k}$.
\end{proof}     

And now the multiplicity result can be proved.
\begin{thm} Let $(M,g)$ be a surface of revolution which is diffeomorphic to $S^{2}$ and let $\lambda_{m}$ be the $m$th distinct eigenvalue of $-\Delta$.
 Then $dimE_{\lambda_{m}}\leq2m+1$ for all $m$, and equality holds for all $m$ if and only if $(M,g)$ is isometric to a sphere of constant curvature.
\end{thm}

\begin{proof}
The second part of the theorem is just Proposition 2.1, iv. above. So we need only prove $dimE_{\lambda_{m}} \leq 2m+1$ for all $m$.

If $m=0$ then $\lambda_{0}=0$ and $dimE_{\lambda_{0}}=dimH^{0}(M,\R)=1$ by Hodge theory.
 So we may assume, henceforth, that $m>0$. From Lemma 3.1 with $j=0$, we have:
\begin{equation} 
\lambda_{m} \leq \lambda_{m}^{1} \; \; (\forall m).    \label{eq:3}
\end{equation} 
Using Proposition 2.1, i. we see that $\lambda_{m} \in SpecL_{k}$ for some $k\in \Z$, and hence for some $k\geq 0$, due to the symmetry exhibited by equation \eqref{eq:1}.
 Finding $dimE_{\lambda_{m}}$ is just a matter of counting the summands in equation \eqref{eq:1}.
 This is done in the following manner. Define $l_{m}=card \left\{ k \in \N | \lambda_{m} \in  SpecL_{k} \right\} $. Then:

$$dimE_{\lambda_{m}}= \left\{ \begin{array}{ll}
 2l_{m} & \mbox{if $\lambda_{m} \not\in SpecL_{0}$} \\
 2l_{m}+1 & \mbox{if $\lambda_{m} \in SpecL_{0}$}
                      \end{array} \right. .$$
So that
\begin{equation}
dimE_{\lambda_{m}} \leq 2l_{m}+1.   \label{eq:4} 
\end{equation}

We can now prove the inequality by contradiction. Suppose, for some $m>0$, that $dimE_{\lambda_{m}} > 2m+1$.
 Then from \eqref{eq:4}, $2m+1<2l_{m}+1$ i.e. $m<l_{m}$. In other words, there are more than $m$ distinct natural numbers in the set 
$\left\{ k \in \N | \lambda_{m} \in SpecL_{k} \right\}$.
 Hence, there exists an element of this set, say $k^{\ast}$, which satisfies $m<k^{\ast}$ while at the same time $\lambda_{m} \in SpecL_{k^{\ast}}$.
 To summarize: $\lambda_{m}\leq \lambda_{m}^{1}$ by \eqref{eq:3}, $\lambda_{m}^{1}< \lambda_{k^{\ast}}^{1}$ because of Proposition 2.1 ii, and $\lambda_{k^{\ast}}^{1} \leq \lambda_{m}$ because $\lambda_{m} \in SpecL_{k^{\ast}}$ and $\lambda_{k^{\ast}}^{1}$ is the first eigenvalue of $L_{k^{\ast}}$.
Putting these three inequalities together produces the contradiction $\lambda_{m}<\lambda_{m}$, and the proof is finished.\end{proof} 

\section{Sharp Bounds for the Eigenvalues}

In the case of Surfaces of Revolution diffeomorphic to $S^2$, sharp upper bounds for the distinct eigenvalues can be found.   

\begin{thm} Let $(M,g)$ be a surface of revolution which is diffeormorphic to $
S^2$ and whose metric is given by \eqref{eq:met} and \eqref{eq:f}. Let $\lambda_{m}$ be the $m$th
 distinct eigenvalue of $-\Delta$, then

\begin{equation}
\lambda_{m}\leq m^2\left[\frac{\int_{-1}^{1}f^{m-1}(x)dx}{\int_{-1}^{1}f^{m}(x)d
x} \right]+\frac{m\int_{-1}^{1}f^{m}(x)K(x)dx}{2\int_{-1}^{1}f^{m}(x)dx}
\label{eq:sharp}
\end{equation}
and equality holds for all $m$ if and only if $(M,g)$ is isometric with $(S^2,can)$.
\end{thm}

\begin{proof} 
From Lemma 3.1 with $j=0$, $\lambda_{m}\leq \lambda_{m}^1 \;\; \forall m$.
 Since $\lambda_{m}^1$ is the first eigenvalue of the operator \eqref{eq:odo}, the minimum principle shows that
$$\lambda_{m}\leq \frac{\int_{-1}^{1}\left[f(x)\left(\frac{du}{dx}\right)^2+
\frac{m^2}{f(x)}u^2\right]dx}{\int_{-1}^{1}u^2dx}$$
$\forall u\in C^{\infty}(-1,1)$ such that $u(-1)=u(1)=0$.
Setting $u=f^{l/2}(x)$, integrating by parts, and using the fact that $K(x)=(-1/2)f''(x)$ yields

\begin{equation}
\lambda_{m}\leq m^2\left[\frac{\int_{-1}^{1}f^{l-1}(x)dx}{\int_{-1}^{1}f^{l}(x)dx}
\right]+\frac{l\int_{-1}^{1}f^{l}(x)K(x)dx}{2\int_{-1}^{1}f^{l}(x)dx} \;\;\forall l\in \N .
 \label{eq:ray}
\end{equation}

The inequality \eqref{eq:sharp} follows immediately by setting $l=m$ in this formula. 

If $(M,g)$ is isometric with $(S^2,can)$ then $f(x)=1-x^2$ and $K(x) \equiv 1$ and a simple calculation shows that
$$\frac{\int_{-1}^{1}(1-x^2)^{m-1}dx}{\int_{-1}^{1}(1-x^2)^{m}dx}=1+\frac{1}{2m}.$$
so that $\lambda_m \leq m^2+m$. But these upper bounds are the eigenvalues for $(S^2,can)$.

Conversely, if equality holds in formula \eqref{eq:sharp} for all $m$, then $f^{m/2}(x)$ is an eigenfunction with eigenvalue $\lambda_m$ and $f^{m}(x)$ is an eigenfunction with eigenvalue $\lambda_{2m}$ so by Proposition 2.1 v., $(M,g)$ is isometric to $(S^2,can)$.
\end{proof}

The following theorem produces explicit sharp upper bounds for a subsequence of the distinct eigenvalues.
  
\begin{thm} Let $(M,g)$ be a surface of revolution of area $4\pi$ which is diffeomorphic to $S^{2}$, and let $\lambda_{m}$ be the $m$th distinct eigenvalue of $-\Delta$.
 For every $k\in \N$ there exists an $m\geq k$ such that $\lambda_{m}\leq m^2+m$. (i.e. there exists a subsequence, $\{m_k\}_{k=1}^{\infty} \subset \N$, such that $\lambda_{m_k}(g)\leq \lambda_{m_k}(can)$).
\end{thm}

\begin{proof} We will show that for every $k\in\N$, there exists $j_k\geq 0$ such that 
\begin{equation}
\lambda_{k}^{j_k +1}\leq (k+j_k)^2 +(k+j_k),   \label{eq:5}
\end{equation}
then, by lemma 3.1, since $\lambda_{k+j_k}\leq \lambda_{k}^{j_k +1}$, the theorem follows by letting $m=k+j_k$.

To prove \eqref{eq:5} we assume the contrary: that for some $k_0$,\\ $\lambda_{k_0}^{j +1} > (k_0 +j)^2 + (k_0 +j)$ for all $j\geq 0$. This means that for all $j$,
$$\frac{1}{\lambda_{k_0}^{j +1}} < \frac{1}{(k_0 +j)^2 + (k_0 +j)}$$
and hence, by the trace formula \eqref{eq:2} 
$$\frac{1}{k_0}= \sum_{j=0}^{\infty} 1/\lambda_{k_0}^{j+1}< \sum_{j=0}^{\infty} 1/[(k_0 +j)^2 +(k_0 +j)]=\frac{1}{k_0},$$
producing the contradiction, $\frac{1}{k_0}< \frac{1}{k_0}$.
\end{proof}

\begin{remark} We suspect that Proposition 2.1, iv., Theorem 3.1, and Theorem 4.2 hold for \underline{all} metrics on $S^2$.
 Some justification for this belief can be found in the contrapositive statements of these results.

\vspace*{.1in}

{\em If there exists a Riemannian manifold $(M,g)$, diffeomorphic to $S^{2}$, which satisfies any one of the conditions
 \begin{description}
   \item[i] $dimE_{\lambda_{m}}=2m+1 \;\; \forall m$ but $\lambda_m(g) \neq \lambda_m(can)$ for some $m$, 
   \item[ii] $dimE_{\lambda_{m}}>2m+1$ for some $m$,
   \item[iii] $\lambda_m(g) \leq \lambda_m(can)$ for only a finite number of eigenvalues,
 \end{description}
 then the group of isometries, $\Im (M,g)$, is finite.}

\vspace*{.1in} 
In item ii. $m \neq 1$ because, as mentioned in the introduction, Cheng \cite{ch1} has proved $dimE_{\lambda_{1}} \leq 3$ for all Riemannian surfaces which are homeomorphic to $S^{2}$.
 Generically, one would not expect to find metrics satisfying any of these conditions. 
There are examples of metrics with ``large" multiplicities, but
 these examples occur in higher dimensions (\cite{bb}, \cite{cdv}, and \cite{ur}) or for surfaces of larger genus (\cite{bu}, \cite{col}).
\end{remark}
\section{Examples}

Although one would expect the best estimates to come from the inequality \eqref{eq:sharp}, the inequality \eqref{eq:ray} with $l=1$ provides rough upper bounds for all the eigenvalues which are much easier to compute and still quite useful.
 They are:

\begin{equation}
\lambda_{m}\leq m^2\left[\frac{2}{\int_{-1}^{1}f(x)d
x} \right]+\frac{\int_{-1}^{1}f(x)K(x)dx}{2\int_{-1}^{1}f(x)dx}.
\label{eq:rough}
\end{equation}

The form of the coefficient of $m^2$ suggests that one might distinguish between the two cases $\int_{-1}^{1}f(x)dx <2$ and $\int_{-1}^{1}f(x)dx \geq 2$.
 In the former case it is easy to see that these rough bounds exceed the eigenvalues of the standard sphere and so provide little new information about the nature of eigenvalues.
 The second case is more interesting as we can see from the following:

\begin{prop} Let $(M,g)$ be a surface of revolution which is diffeormorphic to 
$S^2$ and whose metric is given by \eqref{eq:met} and \eqref{eq:f}, and let $\lambda_{m}$ be the $m$th
 distinct eigenvalue of $-\Delta$. If $\int_{-1}^{1}f(x)dx \geq 2$ then $K(p)<0$ for some $p \in M$ and  
\begin{equation}
\lambda_{m}\leq m^2
+\frac{\int_{-1}^{1}f(x)K(x)dx}{2\int_{-1}^{1}f(x)dx}.
\label{eq:sp}
\end{equation}

Consequently,
$$\lambda_m(g) < \lambda_m(can) \; \; \mbox{ for $m$ sufficiently large,} $$
and
$$\lim_{m \rightarrow \infty} \left( \lambda_m(can)-\lambda_m(g) \right) = \infty.$$
\end{prop}

\begin{proof} Integrating $\int_{-1}^{1}f(x)dx$ by parts twice leads to the identity:
$$\int_{-1}^{1}f(x)dx= 2-\int_{-1}^{1}x^2K(x)dx.$$
So if $\int_{-1}^{1}f(x)dx \geq 2$ then $\int_{-1}^{1}x^2K(x)dx \leq 0$ and hence $K(x_0)<0$ for some $x_0 \in [-1,1]$.

The inequalities and the limit formula follow immediately from \eqref{eq:rough}.
\end{proof}

Although these methods do not prove it Proposition 5.1 suggests that metrics on $S^2$ whose curvature has variable sign have spectra which diverge away from that of the standard sphere.

It is easy to find examples for which \underline{all} of the eigenvalues are less than, and diverging from, those of the standard sphere.

\begin{definition} The function ${\displaystyle f(x)=\frac{2(1-x^2)}{1+x^2}}$ satisfies the conditions \eqref{eq:f} and so defines a surface of revolution metric on $S^2$ of area $4\pi$. It is an elementary exercise to verify that
$$ \int_{-1}^{1}f(x)dx = 2\pi-4, $$
\begin{equation}
K(x)=4\frac{1-3x^2}{(1+x^2)^3} \; \; \mbox{and} \label{eq:Kurv}
\end{equation}
$$ \int_{-1}^{1}f(x)K(x)dx = \pi+\frac{4}{3}.$$
From \eqref{eq:Kurv} we see that this metric has negative curvature on the polar regions defined by the union of intervals $(-1,-1/\sqrt{3})\cup (1/\sqrt{3},1)$ and from \eqref{eq:rough} we have 
$$ 
\lambda_{m}\leq \frac{1}{\pi-2}m^2
+\frac{3\pi+4}{12\pi-24}.$$
As a result
$$\lambda_{m} < m^2+1, \;\;  \forall m.$$
\end{definition}

 \end{document}